\documentclass[nocap]{ctexart}
\usepackage{amsfonts}
\begin{document}
\title{\textsc{Note on Some Limit Properties of Trigonometric Series of A Certain Type}}\author{Yin Li\\Depertment of Mathematics, NanJing Normal University\\210046, China\\\emph{Email address}: \texttt{yinlee1004@sina.com}}\date{}\maketitle
\begin{abstract}
In this note, we study a certain class of trigonometric series which is important in many problems. An unproved statement in Zygmund's book [5] will be proved and generalized. Further discussions based on this problem will also be made here.
\end{abstract}
\begin{flushleft}
\section{Introduction}
For convenience, we identify the torus $\mathbb{T}=\mathbb{R}/2\pi\mathbb{Z}$ with the interval $[-\frac{1}{2},\frac{1}{2}]$. In this note, attention will be paid to the following class of trigonometric series\\
\begin{equation}\sum_{-\infty}^{\infty}a_{n}e^{int},t\in\mathbb{T},\end{equation}\\
where $\{a_n\}$ is a convex even sequence of positive numbers which satisfies\\
\begin{equation}a_{n}\log n=O(1).\end{equation}
Let\\
\begin{equation}f(t)=\sum_{j=0}^{\infty}(j+1)(a_{j}+a_{j+2}-2a_{j+1})F_{j}(t),\end{equation}\\
where $F_{j}(t)$ is the $j$th Fejer kernel. It is not hard to verify that (1) defines the Fourier series of $f(t)$ and $f(t)\in L^{1}(\mathbb{T})$.\\
By the logarithmical growth of $\|D_{N}(t)\|_{L^1(\mathbb{T})}$, where $D_{N}(t)$ denotes the $N$th Dirichlet kernel, one is able to deduce the following proposition ([5], Chapter V, Theorem (1.12)):\\
\textbf{Proposition 1.} \textit{Let} $\{c_{n}\}$ \textit{be a convex even sequence of positive numbers, then the partial sums of the trigonometric series} $\sum_{-\infty}^{\infty}c_{n}e^{int}$ \textit{are bounded in} $L^{1}(\mathbb{T})$ \textit{if and only if it satisfies} (2) \textit{and the series} $\sum_{-\infty}^{\infty}c_{n}e^{int}$ \textit{converges in} $L^{1}(\mathbb{T})$ \textit{if and only if}\\
\begin{equation}c_{n}\log n=o(1).\end{equation}\\
Let $S_{N}(f,t)$ be the $N$th partial sum of the Fourier series of $f$ at $t\in\mathbb{T}$. By Proposition 1, we have\\
\begin{equation}\int_{\mathbb{T}}|S_{N}(f,t)|dt\leq C_1,\end{equation}\\
\begin{equation}\int_{\mathbb{T}}|f(t)-S_{N}(f,t)|dt\leq C_2,\end{equation}\\
where $C_1$ and $C_2$ are some positive constants.\\
Zygmund, in [5], page 185, asserts without proof that instead of (5) and (6), we actually have both $\lim_{N\rightarrow\infty}\int_{\mathbb{T}}|S_{N}(f,t)|dt$ and $\lim_{N\rightarrow\infty}\int_{\mathbb{T}}|f(t)-S_{N}(f,t)|dt$ exist when $f(t)$ has $\sum_{n=2}^{\infty}\frac{\cos nt}{\log n}$ as its Fourier series.\\
It is the purpose of this note to prove a more general fact which contains this assertion as its special case. Further discussions on an interesting problem induced by this assertion will be made in Section 3.\\
\section{Proof of the Generalized Assertion}
We will show the following theorem is true, therefore Proposition 1 can be generalized.\\
\textbf{Theorem 1.} \textit{For every} $f$ \textit{defined by} (3), \textit{both} $\lim_{N\rightarrow\infty}\int_{\mathbb{T}}|S_{N}(f,t)|dt$ \textit{and} $\lim_{N\rightarrow\infty}\int_{\mathbb{T}}|f(t)-S_{N}(f,t)|dt$ \textit{exist}.\\
We begin with a useful property of the Dirichlet kernel.\\
\textbf{Lemma 1.} \textit{Let} $0\leq t_{1}^{2}<t_{2}^{2}<...<t_{N+1}^{2}\leq\frac{1}{2}$ \textit{denote the} $N+1$ \textit{extrema of} $D_{N}(t)$ \textit{on the interval} $[0,\frac{1}{2}]$, \textit{then there exist positive constants} $c_{1},c_{2},...,c_{N+1}$ \textit{such that}\\
\begin{equation}D_{N}(t_{i}^{2})=c_{i}N\textrm{ \textit{and} }\sum_{i=1}^{N+1}c_{i}=O(\log N).\end{equation}\\
\textit{Proof}. Let $0\leq t_{1}^{1}<t_{2}^{1}<...<t_{N+1}^{1}\leq\frac{1}{2}$ denote the points of the set $\{t\in [0,\frac{1}{2}]|D_{N}(t)=\pm\frac{1}{\sin \pi t}\}$. It is obvious that $t_{i}^{2}<t_{i}^{1}<t_{i+1}^{2}$ for $i=1,2,...,N$. Since $\frac{1}{\sin \pi t}$ is increasing and $-\frac{1}{\sin \pi t}$ is decreasing on $[0,\frac{1}{2}]$, we have\\ \begin{equation}|D_{N}(t_{i+1}^{2})|<|D_{N}(t_{i}^{1})|<|D_{N}(t_{i}^{2})|,i=1,2,...,N.\end{equation}\\
We compute\\
\begin{equation}|D_{N}(t_{i}^{1})|=\frac{2}{\sin\frac{2i-1}{4N+2}\pi}.\end{equation}
(7) follows from (8), (9) and the well-known facts $\sin\frac{1}{N}=O(\frac{1}{N})$ and $\sum_{i=1}^{N}\frac{1}{i}=O(\log N)$. $\square$\\
\textsc{Remark}. Similarly, we can define $c_{-N},c_{-N+1},...,c_0=c_1$ for extrema of $D_N(t)$ on the interval $[-\frac{1}{2},0]$.\\

\textbf{Proof of Theorem 1.} Summation by parts twice yields\\
\begin{eqnarray}f(t)-S_{N}(f,t)=\sum_{j=N-1}^{\infty}(j+1)(a_{j}-2a_{j+1}+a_{j+2})F_{j}(t)\\
-NF_{N-1}(t)(a_{N-1}-a_{N-2})+D_{N}(t)a_N.\nonumber
\end{eqnarray}\\
Since the first and the second term in (10) tend to 0 as $N\rightarrow\infty$, we get\\
\begin{equation}\int_{\mathbb{T}}|f(t)-S_{N}(f,t)|dt=a_{n}O(\log n)+o(1).\end{equation}\\
By our assumption (2), this proves $\lim_{N\rightarrow\infty}\int_{\mathbb{T}}|f(t)-S_{N}(f,t)|dt$ exists.\\
To show that $\lim_{N\rightarrow\infty}\int_{\mathbb{T}}|S_{N}(f,t)|dt$ exists, we estimate $\|S_{N}(f,t)\|_{L^{1}(\mathbb{T})}$. We do the partition\\ \begin{equation}\mathbb{T}=\bigcup_{i=-N}^{N-1}[\frac{i}{2N+1},\frac{i+1}{2N+1})\bigcup[-\frac{1}{2},-\frac{N}{2N+1})\bigcup[\frac{N}{2N+1},\frac{1}{2}]=\bigcup_{j=1}^{2N+2}I_j\end{equation}\\ and estimate the order of every integral $\int_{I_j}|S_{N}(f,t)|dt$. By Lemma 1, $\forall\varepsilon>0$, if $N$ is large enough, there is a positive constant $K$ independent of $j$ such that\\
\begin{equation}(K-\varepsilon)|I_j|c_{\sigma(j)}(\sum_{n=2}^{N}\frac{1}{\log n})\leq\int_{I_j}|S_{N}(f,t)|dt\leq(K+\varepsilon)|I_j|c_{\sigma(j)}(\sum_{n=2}^{N}\frac{1}{\log n}),\end{equation}\\
where $c_{\sigma(j)}$ are defined in Lemma 1 and $\sigma$ is a permutation of $\{-N,-N+1,...,N+1\}$ such that $D_N(t_j)=c_{\sigma(j)}N$, $t_j$ being the central point of $I_j$. Thus we have $\sum_{j=1}^{2N+2}c_{\sigma(j)}=O(\log N)$.\\
Since $\sum_{i=2}^{N}\frac{1}{\log i}=O(\frac{N}{\log N})$, sum over the intervals $I_j$ to get\\
\begin{equation}\sum_{j=1}^{2N+2}\int_{I_j}|S_{N}(f,t)|dt=O(\log N)\cdot O(\frac{1}{N})\cdot O(\frac{N}{\log N})=O(1).\end{equation}\\
The proof is complete. $\square$\\

\section{Existence of the Exceptional Set}
By Theorem 1, it is natural to ask whether\\
\begin{equation}\lim_{N\rightarrow\infty}\int_{E}|S_{N}(f,t)|dt\end{equation}\\
exists for every measurable subset $E$ of $\mathbb{T}$.\\
We say that $E\subset\mathbb{T}$ is an \textit{exceptional set} if (15) does not exist. In this section, we shall give an existential proof for the existence of an exceptional set.\\
For a measurable set $E\subset\mathbb{T}$, we have the representation\\
\begin{equation}E=(\bigcup_{j\in\mathbb{Z}_{+}}I_{j})\bigcup N_{1}\bigcap N_{2},\end{equation}\\
where $I_i\bigcap I_j=\emptyset$ provided that $i\neq j$, each $I_j$ is a closed interval or an empty set and $N_1$, $N_2$ are null sets. Thus we shall identify $E$ with $\bigcup_{j\in\mathbb{Z}_{+}}I_{j}$ when considering the existence of (15).\\
The existence of an exceptional set is not evident at first glance. Note that if $E\subset\mathbb{T}$ is a closed interval containing 0, the method we used in estimating $\|S_{N}(f,t)\|_{L^{1}(\mathbb{T})}$ in the proof of Theorem 1 can be applied here to show the existence of (15). Thus (15) also exists when $E\subset\mathbb{T}\setminus\{0\}$ is a closed interval. Hence (15) exists for every interval of $\mathbb{T}$. The following result is a generalization of these observations.\\
\textbf{Theorem 2.} \textit{For a measurable set} $E\subset\mathbb{T}$, (15) \textit{exists in the following two cases}:\\
(\textit{i}) \textit{There exists a} $j\in\mathbb{Z}_{+}$ \textit{such that} $0\in I_j$.\\
(\textit{ii}) $0\in(\mathbb{T}\setminus\bigcup_{j\in\mathbb{Z}_{+}}I_{j})^{o}$, \textit{where} $(\mathbb{T}\setminus\bigcup_{j\in\mathbb{Z}_{+}}I_{j})^{o}$ \textit{denotes the interior of} $\mathbb{T}\setminus\bigcup_{j\in\mathbb{Z}_{+}}I_{j}$.\\
The easiest way to prove Theorem 2 is to use the fact that the Fourier series defined by (1) converges everwhere to $f(t)$ on $\mathbb{T}\setminus\{0\}$. See [5], Chapter V, Theorem (1.5). However, we shall give a different proof of this result here, which is much easier than that of [5].\\
\textbf{Lemma 2.} \textit{Let} $f$ \textit{be as in} (3), \textit{then} $\{S_{N}(f,t)\}$ \textit{converges everywhere to} $f(t)$ \textit{on} $\mathbb{T}\setminus\{0\}$.\\
\textit{Proof}. Suppose $I\subset\mathbb{T}\setminus\{0\}$ is a closed interval. Since $(\frac{\sin(N+1)\pi t}{\sin \pi t})^{2}$ is bounded on $[\delta,\frac{1}{2}]$ for every $0<\delta\leq\frac{1}{2}$, there exists a positive constant $C_3$ such that\\
\begin{equation}|F_j(t)|\leq\frac{C_3}{j+1},j\in\mathbb{Z}_{+},t\in I.\end{equation}\\
Since the factor $a_j+a_{j+2}-2a_{j+1}$ is bounded by $C_{4}\frac{1}{j(\log j)^{2}}$ for $j\geq 2$, where $C_4$ is a positive constant, and the series $\sum_{j=2}^{\infty}\frac{1}{j(\log j)^{2}}$ converges, it follows that uniform convergence of the right side of (3) holds on $I$, thus $f$ is continuous on $I$.\\
Since $\{D_{N}(t)\}$ is uniformly bounded on $I$ and $\{a_n\}$ is of bounded variation, it follows that $\{S_{N}(f,t)\}$ converges uniformly on $I$. By the continuity of $f$ on $I$, $\{S_{N}(f,t)\}$ converges to $f(t)$ for every $t\in I$ (See [2], Proposition 3.3.2).\\
Since for every $t_{0}\in\mathbb{T}\setminus\{0\}$, we can choose a closed interval $I_0\subset\mathbb{T}\setminus\{0\}$ such that $t_{0}\in I_0$, the proof is complete. $\square$\\
\textbf{Proof of Theorem 2.} Since $\{S_{N}(f,t)\}$ converges uniformly to $f(t)$ on $\bigcup_{i\neq j}I_{i}$, it follows that $\lim_{N\rightarrow\infty}\int_{\bigcup_{i\neq j}I_{i}}|S_{N}(f,t)|dt$ exists. Since we have already known that $\lim_{N\rightarrow\infty}\int_{I_j}|S_{N}(f,t)|dt$ exists, (i) follows.\\
Since $\{S_{N}(f,t)\}$ converges uniformly to $f(t)$ on $\bigcup_{j\in\mathbb{Z}_{+}}I_{j}$, (ii) follows. $\square$\\
Now let's turn to the proof of the existence of an exceptional set. The following lemma is critical since it associates our problem with the concept of uniform integrability. We'll say that a family of functions $\mathcal{F}\subset L^{1}(\mathbb{T})$ is uniformly integrable if $\mathcal{F}$ has uniformly absolutely continuous integrals, this is justified by Proposition 4.5.3 in [1].\\
\textbf{Lemma 3.} \textit{Let} $(X,\mathfrak{M},\mu)$ \textit{be a positive measure space}. \textit{If} $\mu(X)<\infty$, $f_n\in L^{1}(X,\mu)$ \textit{and} $\lim_{n\rightarrow\infty}\int_{E}f_nd\mu$ \textit{exists for every} $E\in\mathfrak{M}$, \textit{then} $\{f_n\}$ \textit{is uniformly integrable}.\\
\textit{Proof}. Define $\rho(A,B)=\int_{X}|\chi_A-\chi_B|d\mu$, where $\chi_A$ is the characteristic function of $A\in\mathfrak{M}$, then $(\mathfrak{M},\rho)$ is a complete metric space. For each $n$ we have\\
\begin{eqnarray}|\int_{A}f_nd\mu-\int_{B}f_nd\mu|=|\int_{X}f_n(\chi_A-\chi_B)d\mu|\\
\leq\int_{X}|f_n||\chi_A-\chi_B|d\mu\nonumber\\
=\int_{X}|f_n|\chi_{A-B}d\mu+\int_{X}|f_n|\chi_{B-A}d\mu\nonumber\\
=\int_{A-B}|f_n|d\mu+\int_{B-A}|f_n|d\mu.\nonumber
\end{eqnarray}\\
Since a single function $f_n\in L^1(X,\mu)$ is uniformly integrable, it follows that $\forall\varepsilon>0$, there exists a $\delta>0$ such that if $\rho(A,B)=\mu(A-B)+\mu(B-A)<\delta$, then\\
\begin{equation}|\int_{A}f_nd\mu-\int_{B}f_nd\mu|\leq\int_{A-B}|f_n|d\mu+\int_{B-A}|f_n|d\mu<\varepsilon.\end{equation}\\
Therefore, the mapping $\phi:M\rightarrow\mathbb{R},\phi(E)=\int_{E}f_nd\mu$ is continuous for every $n$.\\
$\forall\varepsilon>0$, for a fixed $N\in\mathbb{Z}_{+}$, let\\
\begin{equation}\mathcal{A}_N=\{E\in\mathfrak{M}:|\int_E(f_n(x)-f_N(x))d\mu|<\varepsilon,n>N\}.\end{equation}\\
Since $\lim_{n\rightarrow\infty}\int_{E}f_nd\mu$ exists for every $E\in\mathfrak{M}$ by our assumption, we have $\mathfrak{M}=\bigcup_{N\in\mathbb{Z}_{+}}\mathcal{A}_N$. Hence the Baire category theorem implies that there exists an $N\in\mathbb{Z}_{+}$ such that $\mathcal{A}_N$ has nonempty interior. This means that $\forall\varepsilon>0$, $E_0\in\mathcal{A}_N$, there exist $\delta>0$, $N\in\mathbb{Z}_{+}$ such that if $\rho(E,E_0)<\delta$, $n>N$ then\\
\begin{equation}|\int_{E}(f_n-f_N)d\mu|<\varepsilon.\end{equation}\\
Therefore, if $\mu(A)<\delta$, we have $|\int_{E_0-A}(f_n-f_N)d\mu|<\varepsilon$ and $|\int_{E_0\bigcup A}(f_n-f_N)d\mu|<\varepsilon$. This implies that\\
\begin{equation}|\int_{A}(f_n-f_N)d\mu|<2\varepsilon.\end{equation}\\
Since that family $\{f_1,f_2,...,f_N\}$ is clearly uniformly integrable, there exists a $\delta'>0$ such that when $\mu(A)<\delta'$ we have\\
\begin{equation}|\int_{A}f_nd\mu|<3\varepsilon\end{equation}\\
for all $n\in\mathbb{Z}_{+}$. This proves the lemma. $\square$\\
\textsc{Remark}. The lemma proved above originates from an exercise in [4].\\
Now we are able to establish the existence of the exceptional set.\\
\textbf{Theorem 3.} \textit{There exists a measurable subset} $E\subset\mathbb{T}$ \textit{such that} (15) \textit{does not exist}.\\
\textit{Proof}. If the exceptional set does not exist, by Lemma 3, the family $\{|S_N(f,t)|\}$ is uniformly integrable, thus the family $\{S_N(f,t)\}$ is uniformly integrable. Since $f(t)\in L^1(\mathbb{T})$, Lemma 2 and the Vitali convergence theorem can be applied to show  that $S_N(f,t)$ converges to $f(t)$ in $L^1$-norm, but this contradicts Proposition 1. $\square$\\
\textsc{Remark}. It is easy to see from Lemma 3 that to prove Theorem 3, we only need to show the family $\{S_N(f,t)\}$ is not uniformly integrable, a fact which is weaker than Lemma 2. Actually, we can show that uniform integrability fails for $\{S_N(f,t)\}$ without using Lemma 2 and the Vitali convergence theorem, this serves as a second proof of Theorem 3.\\
To see this, fix an integer $N_0>0$. We estimate the integral $\int_{Q}S_n(f,t)dt$ with $Q\subset\mathbb{T}$ and $m(Q)=\frac{2}{2N_0+1}$, where $m$ represents the normalized Haar measure on $\mathbb{T}$.\\
For an arbitrary $n\in\mathbb{Z}_{+}$, we may assume $n\in(bN_0,(b+1)N_0]$ with a unique $b\in\mathbb{Z}_{+}$. Partition $\mathbb{T}$ as in (12), i.e., $\mathbb{T}=\bigcup_{j=1}^{2n+2}I_j$, and use $J_1,J_2,...,J_{n+2}$ to denote the intervals on which $D_n(t)$ is nonnegative. For $i=1,2,...,n+2$, define\\
\begin{equation}d_i=d(0,J_i)=\sup_{t\in J_i}|t|.\end{equation}\\
We may assume without loss of generality that $d_1\leq d_2\leq...\leq d_{n+2}$. Take $Q=\bigcup_{i=1}^{b-1}J_i\bigcup J$, where $J$ is a measurable subset of $J_b$. For $N_0$ large enough we have\\
\begin{eqnarray}|\int_{Q}S_n(f,t)dt|=\sum_{i=1}^{b-1}\int_{J_i}S_n(f,t)dt+C_5\\
\geq C_6\sum_{i=1}^{b-1}\frac{2}{2bN_0+1}\cdot(c_{2i-1}\sum_{j=2}^{bN_0}\frac{1}{\log j})+C_5\nonumber\\
=\frac{2C_6}{2bN_0+1}(\sum_{i=1}^{b-1}c_{2i-1})\cdot(\sum_{j=2}^{bN_0}\frac{1}{\log j})+C_5\nonumber\\
=O(\frac{1}{bN_0})\cdot O(\log b)\cdot O(\frac{bN_0}{\log bN_0}),\nonumber
\end{eqnarray}\\
where $C_5$ and $C_6$ are positive constants, $c_i$ are as in Lemma 1.\\
Take $b=N_0$, we get $|\int_{Q}S_n(f,t)dt|\geq O(1)$. Thus the family $\{S_n(f,t)\}$ is not uniformly integrable.\\

\section{Remarks}
\textbf{1.} Let $\mathfrak{M}$ be the $\sigma$-algebra formed by the measurable subsets of $\mathbb{T}$, and let $\rho$ be the metric defined in the proof of Lemma 3. Use $\mathfrak{E}$ to denote the family of exceptional sets of (15). Since every interval of $\mathbb{T}$ belongs to the family $\mathfrak{M}\setminus\mathfrak{E}$, we have the following result:\\
\textbf{Proposition 2.} $\mathfrak{E}$ \textit{is of first category in the metric space} $(\mathfrak{M},\rho)$.\\
Proposition 2 shows the reason why it is easier to prove the existence of an exceptional set than to construct one.\\
\textbf{2.} Analogously, one can ask whether $\lim_{N\rightarrow\infty}\int_{E}|f(t)-S_N(f,t)|dt$ exists for every measurable subset $E$ of $\mathbb{T}$. Use an argument similar to that in the proof of Theorem 3, it is easy to see that there exists an exceptional set $E\subset\mathbb{T}$ such that $\lim_{N\rightarrow\infty}\int_{E}|f(t)-S_N(f,t)|dt$ does not exist.\\

\end{flushleft}
\end{document}